\documentclass[12pt]{article}

\usepackage{amssymb,amsmath}
\usepackage{graphicx}
\usepackage{multirow}
\usepackage{geometry}
 \usepackage[pdftex,
 	colorlinks,
 	pdftitle={Higher order variational time discretization of optimal control problems},
 	pdfauthor={C.M. Campos, O. Junge, S. Ober-Blöbaum},
 	hyperindex,
 	filecolor=blue,
 	urlcolor=blue,
 	citecolor=blue,
 	linkcolor=blue]{hyperref}

\newcommand{\bbR}{\mathbb{R}}
\newcommand{\calC}{\mathcal{C}}
\newcommand{\IntAct}{\mathfrak{G}}
\newcommand{\IntFrc}{\mathfrak{F}}
\newcommand{\ObjFun}{\mathfrak{J}}

\newcommand{\diff}{\mathrm{d}}
\newcommand{\dt}{\diff t}
\newcommand{\dd}[2][{}]{\frac{\diff#1}{\diff#2}}
\newcommand{\ddt}{\dd[]t}
\newcommand{\pp}[2][{}]{\frac{\partial#1}{\partial#2}}
\newcommand{\dq}{\dot{q}}
\newcommand{\cf}{\emph{cf.}}
\newcommand{\eg}{\emph{e.g.}}

\newcommand{\ie}{\emph{i.e.}}

\newcommand{\quand}{\quad\textrm{and}\quad}

\title{
  Higher order variational time discretization\\ of optimal control problems%
  \thanks{
    Funding for this paper furnished by the EU Marie Curie initial training
    network \emph{Sensitivity analysis for deterministic controller design}
    (FP 7) and the DFG Collaborative Research Center 614
    \emph{Self-Optimizing Con\-cepts and Struc\-tures in Mechanical
    Engineering}
  }
}

\date{April 2012}

\author{C.M.~Campos \and O.~Junge\thanks{Faculty for Mathematics, Technische Universit\"at M\"unchen, \texttt{cedricmc@ma.tum.de},\texttt{oj@tum.de}} \and
  S.~Ober-Bl\"obaum\thanks{Department of Mathematics, University of Paderborn, \texttt{sinaob@math.upb.de}}
}

\begin{document}

\baselineskip16pt

\maketitle

\begin{abstract}
We reconsider the variational integration of optimal control problems for mechanical systems based on a direct discretization of the \emph{Lagrange-d'Alembert principle} as proposed in \cite{ObJuMa10}.  This approach yields discrete dynamical constraints which by construction preserve important structural properties of the system, like the evolution of the momentum maps or the energy behaviour.

Here, we employ higher order quadrature rules based on polynomial collocation.  The resulting variational time discretization decreases the overall computational effort.
\end{abstract}

% INTRODUCTION -----------------------------------------------------------------
\section{Introduction}
In recent years, much effort in designing numerical methods for the time integration of (ordinary) differential equations has been put into schemes which are \emph{structure preserving} in the sense that important \emph{qualitative} features of the original dynamics are preserved in its time discretization, \cf\ the recent monograph \cite{HaLuWa02}.  A particularly elegant way to, e.g., derive symplectic integrators is by discretizing Hamilton's principle as suggested by \cite{Su90a,Ve88a}, see also \cite{MaWe01}.

Evidently, structure preservation might equally be important in optimal control problems. In fact, in \cite{ObJuMa10} a new approach\footnote{Referred to as DMOC as in \emph{Discrete Mechanics and Optimal Control}.} to the transcription of a mechanical optimal control problem into a finite dimensional nonlinear programming problem has been proposed which is based on a direct discretization of the \emph{Lagrange-d'Alembert principle} (instead of the associated Euler-Lagrange differential equations of motion).  This approach yields a finite-difference type discretization of the dynamical constraints of the problem which by construction preserves important structural properties of the system, like the evolution of the momentum maps associated to the symmetries of the Lagrangian or the energy behaviour \cite{MaWe01,ObJuMa10}.

So far, quadrature rules of second order have been used in order to approximate the action functional of the system.  In this work, we employ higher order rules based on polynomial collocation as suggested in \cite{MaWe01} for variational integrators. This decreases the overall computational effort of the approach, while maintaining its structure preservation properties.

The paper is structured as follows: In Section~\ref{sec:oc} we introduce the mechanical optimal control problem under consideration, while in Section~\ref{sec:discrete} we describe its variational discretization using higher order polynomial collocation.  We present two numerical experiments in Section~\ref{sec:experiments} and conclude by outlining future research directions in Section~\ref{sec:conclusion}.
% ------------------------------------------------------------------------------

% MECHANICAL OPTIMAL CONTROL PROBLEMS ------------------------------------------
\section{Mechanical optimal control problems}\label{sec:oc}
We consider a mechanical system with \emph{configuration manifold} $Q$ together with a $\calC^k$ \emph{Lagrangian} $L\colon TQ\to\bbR$, $k\geq2$, where the associated \emph{state space} $TQ$ describes the position and velocity of a particle moving in the system. Usually, the Lagrangian takes the form of kinetic minus potential energy, $ L(q,\dq) = K(q,\dq) - V(q) = \frac12\,\dq^T\cdot M(q)\cdot\dq - V(q)$, for some (positive definite) \emph{mass matrix} $M(q)$.

Given a time interval $[0,T]$, the Lagrangian defines an \emph{integral action} $\IntAct\colon\calC^2([0,T],\allowbreak Q)\to\bbR$ on the space of twice differentiable paths $q\colon[0,T]\to Q$, namely:
\[ \IntAct(q) = \int_0^T L(q(t),\dq(t)) \,\dt \,.\]
Note that $\calC^2([0,T],Q)$ is a infinite dimensional smooth manifold and that $\IntAct$ is of class $\calC^k$. To be rigorous, $\calC^2([0,T],Q)$ should be completed with some norm to form a Banach space, although this falls out of the scope of the present work.

The \emph{prinpicle of least action}, also known as \emph{Hamilton's principle}, establishes that a ``particle'' moving in $Q$ from a point $q^0\in Q$ to a point $q^T\in Q$ will do so along a path $q\colon[0,T]\to Q$ such that $q(0)=q^0$ and $q(T)=q^T$ and $q$ itself is an extremal of $\IntAct$. An \emph{extremal} of the integral action $\IntAct$ is a trajectory $q\in\calC^2([0,T],Q)$ so that, for any infinitesimal variation $\delta q\in T\calC^2([0,T],Q)=\calC^1([0,T],TQ)$ of $q$ (\ie\ $\tau_Q\circ\delta q=q$, where $\tau_Q\colon TQ\to Q$ is the canonical projection), we have that
\[ \delta\IntAct(q)\cdot\delta q = 0 \,. \]
Using integration by parts and considering infinital variations that are null at the end points, $\delta q(0)=\delta q(T)=0$, it is easy to show that the previous expression is equivalent to the following one
\[ \pp[L]{q} - \ddt\pp[L]{\dq} = 0 \,, \]
which is the well known \emph{Euler-Lagrange equation}.

The trajectories resulting from the previous second order differential equation are free in the sense that no force exerts any influence on the mechanical system, which is a rather simplistic assumption. Otherwise and more generally, if an external force $f\colon TQ\to T^*Q$ exerts some influence on the system, the principle of least action is then know as the \emph{principle of virtual work} or \emph{Lagrange-d'Alembert principle}, which states that in presence of an external force the extremal trajectories must satisfy the variational equation
\begin{equation} \label{eq:LDA-principle}
\delta \int_0^T L(q(t),\dq(t))\,\dt + \int_0^T f(q(t),\dq(t))\cdot \delta q(t) \,\dt = 0 \,,
\end{equation}
for any infinitesimal variation $\delta q$ that is null at the end points. This equation is in turn equivalent to the \emph{forced Euler-Lagrange equation}
\begin{equation} \label{eq:LDA-equation}
\pp[L]{q} - \ddt\pp[L]{\dq} + f = 0 \,.
\end{equation}

When the Lagrangian is \emph{regular}, that is when the velocity Hessian matrix $\partial^2 L/\partial \dq^2$ is non-degenerate, the Euler-Lagrange equation \eqref{eq:LDA-equation} has a unique solution $q\in\calC^2([0,T],Q)$ for each initial condition $(q^0,\dq^0)\in TQ$ (rather than for boundary conditions $q^0,q^T\in Q$). Therefore, it defines a map $F_L\colon TQ\times\bbR\to TQ$ by $F_L(q^0,\dq^0,t) := (q(t),\dq(t))$. 
In fact, $F_L$ is a $\calC^{k-1}$-diffeomorphism (for fixed $t\in\bbR$) called the \emph{Lagrangian flow} of the system.

Imagine now that we want to control the trajectories of the system following a minimal cost criteria.  We must then first assume that the mechanical system may be driven by means of some \emph{control parameter} $u\in U\subset\bbR^\ell$ on which the external force $f$ will therefore depend. We also assume that an infinitesimal \emph{cost function} $C\colon TQ\times_QT^*Q\to\bbR$ is given associated to the \emph{objective functional} $\ObjFun\colon\calC^2([0,T],Q)\times\calC([0,T],U)\to\bbR$ 
\begin{equation} \label{eq:objective.functional}
\ObjFun(q,u) = \int_0^T C(q(t),\dq(t),f(q(t),\dq(t),u(t))) \,\dt \,.
\end{equation}

We seek for pairs of curves $(q,u)\colon[0,T]\to Q\times U$ such that $q$, under the influence of the external force $f(q(\cdot),\dq(\cdot),u(\cdot))$, moves from a given initial state $(q^0,\dq^0)\in TQ$ to a given final state $(q^T,\dq^T)\in TQ$ while minimizing the objective functional $\ObjFun$. That is, we seek to solve the \emph{mechanical optimal control problem}
\begin{subequations} \label{eq:mocp}
\begin{align}
\label{eq:mocp-min}
& \min_{q,u}\ObjFun(q,u) \,,\ q\in\calC^2([0,T],Q) \,,\ u\in\calC([0,T],U) \\
\label{eq:mocp-dyn}
\textrm{s.t.}\ 
& \pp[L]{q} - \ddt\pp[L]{\dq} + f = 0 \\
\label{eq:mocp-cond.ini}
& (q(0),\dq(0)) = (q^0,\dq^0), \quad  (q(T),\dq(T)) = (q^T,\dq^T) 
\end{align}
\end{subequations}
 
As stated above, if $L$ is regular, then \eqref{eq:mocp-dyn} defines a unique solution $q\in\calC^2([0,T],\allowbreak Q)$ (for fixed $u\in C([0,T],U)$) satisfying the initial condition \eqref{eq:mocp-cond.ini}. Note that the assumption on $u$ to be continuous is rather restrictive from the point of view of typical applications.  Since the focus of this paper is on the approximation of the solution curve $q$ by (higher order) polynomials we nevertheless restrict our attention to this case here. In the following we assume the regularity of $L$ as well as the existence of an optimal solution $(q^*,u^*)$ to the mechanical optimal control problem and focus on the integration of the dynamical part in order to solve the optimal control problem numerically.

% HIGHER ORDER VARIATIONAL DISCRETIZATION --------------------------------------
\section{Higher order variational discretization}\label{sec:discrete}
Solution methods for optimal control problems can be divided into \emph{indirect and direct approaches} (cf.~\cite{Binder01}). While the indirect approach bases on the solution of the necessary optimality conditions, the direct approach transforms the problem into a finite dimensional restricted optimization problem by a discretization of the forced Euler-Lagrange equation \eqref{eq:LDA-equation}. In this work, we follow a direct approach, however instead of dscretizing the equation of motion, we discretize the variational principle \eqref{eq:LDA-principle} (cf.~\cite{ObJuMa10}).

In order to do so, we use piecewise polynomial approximations to the trajectories and numerical quadrature to approximate the integrals following \cite{MaWe01}. To this end, we divide the interval $[0,T]$ into smaller subintervals $I_k=[t_k,t_{k+1}]$ ($t_0=0$, $t_N=T$) of fix length $h=T/N$ and on each of these subintervals we perform the following construction: We approximate $q\colon I_k \to Q$ and $u\colon I_k\to U$ by polynomials $q_k^d\colon[0,h]\to Q$ and $u_k^d\colon [0,h]\to U$ of degree $s$ and $m$, respectively. Given intermediate times $0=d_0<d_1<\dots<d_{s-1}<d_s=1$ and intermediate points $q_k=q_k^0,q_k^1,q_k^2,\dots,q_k^{s-1},q_k^s=q_{k+1}$,  the interpolating polynomial $q_k^d$ of degree $s$ with $q_k^d(d_\nu h)=q_k^\nu$ for $\nu =0,\dots,s$ is uniquely defined. Analogously, we choose a set of interior times $0=\tilde d_0<\tilde d_1<\dots<\tilde d_{m-1}<\tilde d_m=1$ and interior points $u_k^0,\dots,u_k^m\in U$ that represent a parametrization of the space of polynomials $u_k^d\colon[0,h]\to U$ of degree $m$. Obviously, to ensure the proper definition of the polynomials implies the assumption of a linear structure on $Q$, canonical or taken, for instance, from a global chart. 

Note that by identifying the polynomials $q_k^d$ with the set of intermediate points $q_k^0,\dots,q_k^s$ (resp. $u_k^d$ with $u_k^0,\dots,u_k^m$), we are implicitely identifying the set of polynomials of order $s$ from $[0,h]$ to $Q$ with $Q^{s+1}$ (resp. with $U^{m+1}$) and, therefore, its tangent space with $\oplus_{s+1}TQ$ (resp. with $\oplus_{m+1}TU$).

As for the derivation of the continuous time dynamical equations, the control parameter plays no role in the derivation of the discrete time dynamical equations.  In order to simplify notation we therefore ommit the explicit dependence on $u$ and also write $f(c_ih)$ in order to denote $f(q(c_ih),\dot q(c_ih),u(c_ih))$ in the following.

\paragraph{Higher order discrete mechanics}
We approximate the action integral on $[0,h]$ by numerical quadrature,
\begin{equation} \label{eq:lagrangian.multipoint}
L_d (q_k^0,q_k^1,\dots,q_k^s) := h\sum\limits_{i=1}^r b_i L(  q_k^d(c_ih),\dq_k^d(c_ih)) \approx \int_{I_k} L(q,\dq) \,\dt
\end{equation}
where $c_i \in [0,1]$, $i=1,\dots, r$, are the quadrature nodes and $b_i$ the corresponding weights. $L_d$ is called a \emph{multipoint discrete Lagrangian}. The \emph{discrete action sum} over the entire trajectory on $[0,T]$ is then
\[ \IntAct_d(\{ (q_k^0,\dots,q_k^s) \}_{k=0}^{N-1}) := \sum\limits_{k=0}^{N-1} L_d(q_k^0,q_k^1,\dots,q_k^s) \approx \int_0^T L(q,\dq) \,\dt \]

Similarly, we approximate the force integral in Equation \eqref{eq:LDA-principle} on $I_k$. Using the short-hand notation $ f(t) = f(q_k^d(t),\dq_k^d(t))$ we define the \emph{multipoint discrete forces} by 
\begin{equation} \label{eq:force.multipoint}
f_k^\nu = f_k^\nu(q_k^0,\dots,q_k^s) := h\sum\limits_{i=1}^r b_i f(c_i h)\frac{\partial q_k^d(c_i h)}{\partial q_k^\nu} \,.
\end{equation}
Then, the force integral is approximated by
\[
F_d(q_k^0,\dots,q_k^s) \cdot \delta(q_k^0,\dots,q_k^s) := \sum\limits_{\nu=0}^s f_k^\nu\cdot\delta q_k^\nu 
= h \sum\limits_{i=1}^r b_i f(c_ih) \cdot \delta q_k^d(c_ih) \approx \int_{I_k} f\cdot \delta q \,\dt \,,
\]
where $\delta q_k^d(t)$ denotes variations of $q_k^d$ given as $\delta q_k^d(t) = \sum\limits_{\nu=0}^s \frac{\partial q_k^d(t)}{\partial q_k^\nu}\delta q_k^\nu$.
Thus, the discrete force integral over the entire trajectory on $[0,T]$ is
\[
\IntFrc_d(\{ (q_k^0,\dots,q_k^s) \}_{k=0}^{N-1}) = \sum\limits_{k=0}^{N-1} F_d(q_k^0,\dots,q_k^s) \approx \int_0^T f\cdot \delta q \,\dt \,.
\]

Having defined the discrete Lagrangian action $\IntAct_d$ and the discrete action of the force $\IntFrc_d$, we require that the discrete version of the Lagrange-d'Alembert principle holds for variations of $q_k^\nu$. That is, a sequence of points $\{ (q_k^0,\dots,q_k^s) \}_{k=0}^{N-1}$ is an extremal trajectory for the system if, for any variation $\{ \delta(q_k^0,\dots,q_k^s) \}_{k=0}^{N-1}$ with $\delta q_k^0=\delta q_{N-1}^s=0$, we have that
\begin{equation} \label{eq:DLDA-principle}
(\delta\IntAct_d + \IntFrc_d) \cdot \{ \delta(q_k^0,\dots,q_k^s) \}_{k=0}^{N-1} = 0 \,.
\end{equation}
A simple computation shows that the extended set of discrete Euler-Lagrange equations can be derived as
\begin{subequations} \label{eq:DLDA-equation}
\begin{align}
\label{eq:DLDA-equation-main}
& D_{s+1} L_d(q_k^0,\dots,q_k^s) + D_1 L_d(q_{k+1}^0,\dots,q_{k+1}^s) + f_k^s+f_{k+1}^0 = 0 \,, \\
\label{eq:DLDA-equation-scnd}
& D_{\nu+1} L_d(q_k^0,\dots,q_k^s) + f_k^\nu = 0,\quad \nu = 1,\dots, s-1 \,.
\end{align}
\end{subequations}
for $k=0,\ldots,N-1$, where
\[ D_{\nu+1} L_d(q_k^0,\dots,q_k^s) = \frac{\partial L_d}{\partial q^\nu} = h\sum\limits_{i=1}^r b_i \left( \frac{\partial L}{\partial q}\frac{\partial q_k^d(c_ih)}{\partial q_k^\nu} + \frac{\partial L}{\partial \dq}\frac{\partial \dq_k^d(c_ih)}{\partial q_k^\nu}\right) \,. \]

\paragraph{Alternative construction}
The previous construction is a direct derivation of the discrete Euler-Lagrange equations, which can also be derived in a two-step construction. Since this alternative construction is a reinterpretation of the previous one, we should not get into much detail.

In a first step, we approximate the continuous integral action and force integral by the already used quadrature for $(c_i,w_i)$, giving rise to (still) continuous action sum and force $\IntAct'$ and $\IntFrc'$. Then we define a discrete Lagrangian $L_d:Q\times Q\to\bbR$ that, for two points $q_0,q_1\in Q$ ``sufficiently close'' and a time step $h>0$, gives the value
\[ L_d(q_0,q_1) := \IntAct'(q^d) \,, \]
where $q^d$ satisfies the Lagrange-d'Alembert principle for $\IntAct'$ and $\IntFrc'$ restricted to the class of polynomials of degree $s$ joining $q_0$ and $q_1$ in time $h$. With the proper interpretation, one may check that the polynomial $q^d$ is in fact characterized by the equations \eqref{eq:DLDA-equation-scnd}. Note that, even if it is not explicitely stated, $q^d$ and therefore $L_d$ depend on the external force $f$.

In the second step, we define a discrete action sum on $Q^{N+1}$ by
\[ \IntAct_d(q_0,\dots,q_N) = \sum_{k=0}^{N-1}L_d(q_k,q_{k+1}) \]
and discrete forward and backward forces
\[
f_d^+(q_0,q_1)\cdot\delta(q_0,q_1) := \IntFrc'(q^d)\cdot\frac{\partial q^d}{\partial q^1}\delta q^1
\quand
f_d^-(q_0,q_1)\cdot\delta(q_0,q_1) := \IntFrc'(q^d)\cdot\frac{\partial q^d}{\partial q^0}\delta q^0 \,,
\]
where $q^d$ is given as before. Applying now the discrete Lagrange-d'Alembert principle to discrete action sum $\IntAct_d$ and force
\[ \IntFrc_d(q_0,\dots,q_N)\cdot\delta(q_0,\dots,q_N) := \sum_{k=0}^{N-1} \left[ f_d^+(q_k,q_{k+1}) + f_d^-(q_k,q_{k+1}) \right]\cdot\delta(q_k,q_{k+1}) \,, \]
for any variation $\delta(q_0,\dots,q_N)$ such that $\delta q_0=\delta q_N=0$, we obtain
\[ D_2L_d(q_{k-1},q_k) + D_1L_d(q_k,q_{k+1}) + f_d^+(q_{k-1},q_k) + f_d^-(q_k,q_{k+1}) = 0 \,, \]
for any $k=1,\dots,N-1$, which is equivalent to \eqref{eq:DLDA-equation-main} (with the obvious identifications).

The main advantage of this alternative construction is that, as it may be seen from a two-point discrete Lagrangian approach, one can make the most of the usual theory. For instance, we may define the forward and backward Legendre transformations
\begin{align*}
& \mathbb{F}^+(q_0,q_1) := (q_1,  D_2L_d(q_0,q_1) + f_d^+(q_0,q_1)) \in T^*Q \,, \\
& \mathbb{F}^-(q_0,q_1) := (q_0, -D_1L_d(q_0,q_1) - f_d^-(q_0,q_1))  \in T^*Q \,,
\end{align*}
which are useful to implement the initial and final condition of the optimal control problem from a momentum description or to give the momenta associated to the macro nodes $q_0,\dots,q_N$ along the trajectory.

\paragraph{Discrete approximation of the objective functional}
The objective functional $\ObjFun$ is approximated by a discrete objective function $\ObjFun_d$ using the same numerical quadrature rule and the same polynomials $q_k^d$ and $u_k^d$ as above. For each time interval $I_K$, we define the \emph{multipoint discrete cost function}
\begin{align*}
C_d(q_k^0,\dots,q_k^s, u_k^0,\dots,u_k^m)
     := h\sum_{i=0}^r b_iC(q_k^d(c_ih),\dq_k^d(c_ih),u_k^d(c_ih))
\approx \!\!\int_{I_k}\!\! C(q,\dq,u) \,\dt ,%\,,
\end{align*}
which defines the \emph{discrete objective function} over the entire trajectory $[0,T]$
\begin{align*}
\ObjFun_d(\{ (q_k^0,\dots,q_k^s, u_k^0,\dots,u_k^m) \}_{k=0}^{N-1})
     := \sum_{k=0}^{N-1} C_d(q_k^0,\dots,q_k^s, u_k^0,\dots,u_k^m)
\approx \!\!\int_0^T\!\! C(q,\dq,u) \,\dt .%\,.
\end{align*}

The discrete version of the optimal control problem is then to minimize $\ObjFun_d$ subject to the discrete equations \eqref{eq:DLDA-equation} and discretized boundary constraints, namely
\begin{subequations} \label{eq:dmocp}
\begin{align}
&\displaystyle \min_{q,u}\ObjFun_d(\{ (q_k^0,q_k^1,\dots,q_k^s, u_k^0,\dots,u_k^m) \}_{k=0}^{N-1}) \\
\textrm{s.t.}\ 
&\displaystyle D_{s+1} L_d(q_k^0,\dots,q_k^s) + D_1 L_d(q_{k+1}^0,\dots,q_{k+1}^s) + f_k^s+f_{k+1}^0 = 0 \\
&\displaystyle D_{\nu+1} L_d(q_k^0,\dots,q_k^s) + f_k^\nu = 0 \,,\quad \nu = 1,\dots, s-1,\ k=0,\dots,N-1 \\
&\displaystyle (q_0^d(0),\dq_0^d(0)) = (q^0,\dq^0), \quad \displaystyle (q_{N-1}^d(T),\dq_{N-1}^d(T)) = (q^T,\dq^T)
\end{align}
\end{subequations}
This restricted optimization problem can be solved by standard optimization techniques such as \eg\ SQP methods.
% ------------------------------------------------------------------------------

% RESULTS ----------------------------------------------------------------------
\section{Results}\label{sec:experiments}

As proposed in \cite{MaWe01} for uncontrolled systems, we use Lagrange polynomials for the construction of $q_k^d$ and $u_k^d$. For both of them, the collocation points will coincide with the quadrature points of the corresponding Lobatto's quadrature, which will be the same quadrature rule to approximate all the integrals: Action, force and objective function. According to this, quadrature points and weights are given by Table \ref{tab:lobatto} after rescaling from $[-1,1]$ to $[0,1]$, \ie\ $d_\nu=\tilde d_\nu=c_{\nu+1}=(x_{\nu+1}+1)/2$ and $b_i=w_i/2$. Polynomials $q_k^d$ (and similarly $u_k^d$) will then take the form
\[ q_k^d(t) = \sum_{\nu=0}^s P_\nu(t/h)\cdot q_k^\nu \,, \]
where $P_\nu$ are the base of Lagrangian polynomials of degree $s$
\[ P_\nu(t) = \mathop{\prod_{\mu=0,\dots,s}}_{\mu\neq\nu}\frac{t-d_\mu}{d_\nu-d_\mu} \,,\quad \nu=0,\dots,s. \]
We then have that
\[ \frac{\partial q_k^d}{\partial q_k^\nu} = P_\nu(t/h) \quand \dq_k^d(t) = \frac1h\sum_{\nu=0}^s\dot P_\nu(t/h)\cdot q_k^\nu \,. \]

{\small
\begin{table}
\centering
\begin{tabular}{|c|c|c|c|c|c|c|c|c|}
\hline
$n$ & \multicolumn{2}{c|}{3} & \multicolumn{2}{c|}{4} & \multicolumn{2}{c|}{5} & \multicolumn{2}{c|}{6}\\
\hline
 & $x_j$ & $w_j$ & $x_j$ & $w_j$ & $x_j$ & $w_j$ & $x_j$ & $w_j$\\
\hline
\multirow{3}{*}{}
&      0 & $4/3$ & $\pm\frac15\sqrt5$ & $5/6$ &                     0 & $32/45$ & $\pm\sqrt{\frac1{21}(7-2\sqrt7)}$ & $\frac1{30}(14+\sqrt7)$ \\
& $\pm1$ & $1/3$ &             $\pm1$ & $1/6$ & $\pm\frac17\sqrt{21}$ & $49/90$ & $\pm\sqrt{\frac1{21}(7+2\sqrt7)}$ & $\frac1{30}(14-\sqrt7)$ \\
&		 & 		 & 					  &       &         $\pm1$ 		  & $1/10$  &                            $\pm1$ & $1/15$\\\hline
\end{tabular}
\ \\\ \\ \caption{Lobatto's quadrature: points and weights}
\label{tab:lobatto}
\end{table}}

We point out here a slight difference with respect to \cite{MaWe01}: In there the authors propose numerical quadratures for the uncontrolled system with $r=s$, obtaining a variational integrator of order $2s-2$; while the numerical quadrature we use has the same number of nodes as the polynomials, $r=s+1$, obtaining a variational integrator of order $2s$ (cf.~\cite{Sa12}). The numerical scheme for the optimal control problem including the states as well as the adjoint variables coming from the necessary optimality conditions then inherits this property and shows from numerical simulations (see below) a convergence order of $2s$, as one might expect.

The numerical scheme is implemented in Matlab (R2011b) using the built-in function \textsf{fmincon} to solve the non-linear programming problem. We ran several numeric simulations of the controlled harmonic oscillator with Lagrangian $L(q,\dq)=\frac12(\dq^2-5q^2)$, force $f(q,\dq,u)=u$ and cost function $C(q,\dq,u)=u^2$. In these simulations, the system is driven from an initial steady state $(0,0)$ to a final steady state $(1,0)$ in $T=5$ seconds time.

Table \ref{tab:harosc:s2} represent the output obtained from simulations where the trajectories are approximated by polynomials of order $s=2$. The number $N$ of macro nodes ranges from $2^2$ to $2^6$; $\mathrm{err}(q)$, $\mathrm{err}(u)$, $\mathrm{err}(\lambda)$ represent the error committed on the trajectory, controls and adjoint variables, respectively (which are the max. at the coincident macro nodes); $\sharp(\textrm{iter.})$ is the number of iterations done by \textsf{fmincon}; and $\sharp(\textrm{constr.})$ and $\sharp(DL)$ are the number of evaluations of the constraints and the Jacobian of the Lagrangian. Since there is no exact explicit solution, for the error comparisons, we consider the data obtained for $N=2^9$ as the exact one. As expected, we may observe that the error convergence is of order $4$.

\begin{table}
\centering
\begin{tabular}{|l|r|r|r|r|r|}%r|r|}
\hline
$\sharp(\textrm{nodes})$   &    $2^2=4$ &    $2^3=8$ &   $2^4=16$ &   $2^5=32$ &   $2^6=64$ \\
\hline
$\mathrm{err}(q)$          & $1.92\cdot 10^{-1}$ & $2.48\cdot 10^{-3}$ & $1.65\cdot 10^{-4}$ & $1.07\cdot 10^{-5}$ & $5.70\cdot 10^{-7}$ \\
\hline
$\mathrm{err}(u)$          & $6.24\cdot 10^{-1}$ & $1.06\cdot 10^{-2}$ & $7.22\cdot 10^{-4}$ & $4.53\cdot 10^{-5}$ & $8.67\cdot 10^{-6}$ \\
\hline
$\mathrm{err}(\lambda)$    & $7.71\cdot 10^{-1}$ & $2.11\cdot 10^{-2}$ & $1.44\cdot 10^{-3}$ & $9.02\cdot 10^{-5}$ & $5.73\cdot 10^{-6}$ \\
\hline
$\sharp(\textrm{iter.})$   &         12 &         18 &         21 &         30 &         26 \\
\hline
$\sharp(\textrm{constr.})$ &        247 &        665 &      1.474 &      3.931 &      6.735 \\
\hline
$\sharp(DL)$               &     10.386 &     51.894 &    221.148 &  1.155.810 &  3.919.962 \\
\hline
\end{tabular}
\ \\\ \\ \caption{Controlled harmonic oscillator (with $s=2$)}
\label{tab:harosc:s2}
\end{table}

Due to the higher order of the numerical scheme, one may easily surpass the machine precision and, therefore, it is difficult to obtain better numerical results. For instance, this is the case for $s=5$, where the order of the method is expected to be $10$, which surpasses machine precision already for a time-step of length $h=0.01$. Besides, the particularities of the non-linear solver in use play also an important role (in our case, \textsf{fmincon}).

In Table \ref{tab:harosc:comparison}, we compare simulations of lower order with smaller time-step, with simulations of higher order with bigger time-step. From the result, we may observe that by increasing the order of the method, we may decrease the number of needed nodes, while conserving a comparable accuracy. The benefits are clear: Less nodes means less variables, lower memory cost, lower cpu time consumption. For instance, for $s=5$ and $N=2^3$, the number of evaluations of the Lagrangian is two orders of magnitudes lower than for $s=2$ and $N=2^8$.

Note that, since the Jacobian of the constraints is provided by finite differences, the number of constraints evaluations scales linearly with the number of nodes. The use of automatic differentiation \cite{ADOL} for the computation of derivatives would decrease the number of evaluations and improve the performance of the numerical scheme.

\begin{table}
\centering
\begin{tabular}{|l|r|r|}
\cline{2-3}
\multicolumn{1}{l|}{}      &      $s=2$ &      $s=5$ \\
\hline
$\sharp(\textrm{nodes})$   &   $2^5=32$ &    $2^2=4$ \\
\hline
$\mathrm{err}(q)$          & $1.07\cdot 10^{-5}$ & $1.16\cdot 10^{-5}$ \\
\hline
$\mathrm{err}(u)$          & $4.53\cdot 10^{-5}$ & $3.41\cdot 10^{-5}$ \\
\hline
$\mathrm{err}(\lambda)$    & $9.02\cdot 10^{-5}$ & $6.45\cdot 10^{-5}$ \\
\hline
$\sharp(\textrm{iter.})$   &         30 &         27 \\
\hline
$\sharp(\textrm{constr.})$ &      3.931 &      1.204 \\
\hline
$\sharp(DL)$               &  1.155.810 &    187.848 \\
\hline
\end{tabular}
\qquad 
\begin{tabular}{|l|r|r|}
\cline{2-3}
\multicolumn{1}{l|}{}      &      $s=2$ &      $s=5$ \\
\hline
$\sharp(\textrm{nodes})$   &  $2^8=256$ &    $2^3=8$ \\
\hline
$\mathrm{err}(q)$          & $1.94\cdot 10^{-7}$ & $4.66\cdot 10^{-7}$ \\
\hline
$\mathrm{err}(u)$          & $7.31\cdot 10^{-6}$ & $2.00\cdot 10^{-5}$ \\
\hline
$\mathrm{err}(\lambda)$    & $6.64\cdot 10^{-7}$ & $9.51\cdot 10^{-7}$ \\
\hline
$\sharp(\textrm{iter.})$   &         37 &         31 \\
\hline
$\sharp(\textrm{constr.})$ &     38.043 &      2.656 \\
\hline
$\sharp(DL)$               & 87.880.098 &    796.848 \\
\hline
\end{tabular}
\ \\\ \\ \caption{Order comparison: $s=2$ \emph{vs.} $s=5$}
\label{tab:harosc:comparison}
\end{table}
% ------------------------------------------------------------------------------

% CONCLUSIONS ------------------------------------------------------------------
\section{Conclusions}\label{sec:conclusion}
In this work, discrete variational mechanics has been applied to numerically solve optimal control problems for mechanical systems. In extension to \cite{ObJuMa10}, the use of higher order quadrature rules for the discretization leads to higher order optimal control schemes. These schemes are, for a prescribed accuracy of the discrete optimal solution, computationally more efficient than lower order schemes. 
For the future, the convergence rates of the optimal control scheme, that have been observed numerically, have to be proven. Furthermore, based on the discrete variational framework to construct variational schemes of arbitrary order, we plan to derive discretization schemes that adapt the order depending on the state of the system.  
% ------------------------------------------------------------------------------

% BIBLIOGRAPHY -----------------------------------------------------------------
% Bibliographic styles: amsalpha, amsplain, ieeetr, siam
\bibliographystyle{abbrv}
\bibliography{CmpsObrJng12-arXiv}

\end{document}